\documentclass[11pt]{article}
\usepackage{latexsym,amssymb,amsmath,geometry,cite,enumerate,float,verbatim}
\newtheorem{theorem}{Theorem}[section]
\newtheorem{lemma}[theorem]{Lemma}

\usepackage{lineno}
\usepackage{setspace}
\allowdisplaybreaks

\begin{document}

\onehalfspace

\title{A Short Proof for a Lower Bound on the Zero Forcing Number}

\author{M. F\"{u}rst\and D. Rautenbach}

\date{}

\maketitle

\begin{center}
{\small 
Institute of Optimization and Operations Research, Ulm University, Germany\\
\texttt{maximilian.fuerst,dieter.rautenbach@uni-ulm.de}}
\end{center}

\begin{abstract}
We provide a short proof of a conjecture of Davila and Kenter concerning a lower bound on the zero forcing number $Z(G)$ of a graph $G$. 
More specifically, we show that $Z(G)\geq (g-2)(\delta-2)+2$
for every graph $G$ of girth $g$ at least $3$ and minimum degree $\delta$ at least $2$.
\end{abstract}

{\small 

\begin{tabular}{lp{13cm}}
\end{tabular}
}

\section{Introduction}

We consider finite, simple, and undirected graphs and use standard terminology.

For an integer $n$, let $[n]$ denote the set of positive integers at most $n$.
For a graph $G$,
a set $Z$ of vertices of $G$ is a {\it zero forcing set} of $G$ 
if the elements of $V(G)\setminus Z$ have a linear order $u_1,\ldots,u_k$ 
such that, for every $i$ in $[k]$, 
there is some vertex $v_i$ in $Z\cup \{ u_j:j\in [i-1]\}$ such that 
$u_i$ is the only neighbor of $v_i$ outside of $Z\cup \{ u_j:j\in [i-1]\}$;
in particular, 
$N_G[v_i]\setminus (Z\cup N_G[v_1]\cup \cdots \cup N_G[v_{i-1}])=\{ u_i\}$ for $i\in [k]$.
The {\it zero forcing number $Z(G)$} of $G$, 
defined as the minimum order of a zero forcing set of $G$,
was proposed by the AIM Minimum Rank - Special Graphs Work Group \cite{aim,hv} 
as an upper bound on the corank of matrices associated with a given graph.
The same parameter was also considered in connection with 
quantum physics \cite{bg,bm,s} and logic circuits \cite{bghsy}.

In \cite{dk} Davila and Kenter conjectured that 
\begin{eqnarray}
Z(G)&\geq &(g-2)(\delta-2)+2\label{eb3}
\end{eqnarray}
for every graph $G$ of girth $g$ at least $3$ and minimum degree $\delta$ at least $2$.
They observe that, for $g>6$ and sufficiently large $\delta$ in terms of $g$,
the conjectured bound follows by combining results from \cite{aim2} and \cite{cs}.
For $g\leq 6$, it was shown in \cite{gprs,gr},
Davila and Henning \cite{dh} showed it for $7\leq g\leq 10$,
and, eventually, 
Davila, Kalinowski, and Stephen \cite{dks} completed the proof.
The proof in \cite{dks} is rather short itself but relies on \cite{gprs,gr,dh}.
While the cases $g\leq 6$ have rather short proofs,
the proof in \cite{dh} for $7\leq g\leq 10$ extends over more than eleven pages and requires a detailed case analysis.
Therefore, the complete proof of (\ref{eb3}) obtained by combining \cite{gprs,gr,dh,dks} is rather long.

In the present note we propose a considerably shorter and simpler proof.
Our approach only requires a special treatment for the triangle-free case $g=4$ \cite{gprs},
involves a new lower bound on the zero forcing number,
and an application of the Moore bound \cite{ahl}.

\section{Proof of (\ref{eb3})}

Our first result is a natural generalization of the well known fact $Z(G)\geq \delta(G)$ \cite{aim},
where $\delta(G)$ is the minimum degree of a graph $G$.
For a set $X$ of vertices of a graph $G$ of order $n$, let 
$N_G(X)=\left(\bigcup\limits_{u\in X}N_G(u)\right)\setminus X$,
$N_G[X]=X\cup N_G(X)$, and
$\delta_p(G)=\min\left\{ |N_G(X)|:X\subseteq V(G)\mbox{ and }|X|=p\right\}$
for $p\in [n]$.
Note that $\delta_1(G)$ equals $\delta(G)$.

\begin{lemma}\label{lemma1}
If $G$ is a graph of order $n$, then $Z(G)\geq \delta_p(G)$ for every $p\in [n]$.
\end{lemma}
{\it Proof:} Let $Z$ be a zero forcing set of minimum order.
Let $u_1,\ldots,u_k$ and $v_1,\ldots,v_k$ be as in the introduction.
Since, by definition, $\delta_p(G)\leq n-p$, the result is trivial for $p\geq k=n-|Z|$,
and we may assume that $p<k$.
As noted above, we have $N_G[v_i]\setminus (Z\cup N_G[v_1]\cup \cdots \cup N_G[v_{i-1}])=\{ u_i\}$ for $i\in [k]$,
which implies that $X=\{ v_1,\ldots,v_p\}$ is a set of $p$ distinct vertices of $G$.
Furthermore, it implies that $|N_G[X]|\leq |Z|+p$,
and, hence, $\delta_p(G)\leq |N_G(X)|=|N_G[X]|-p\leq |Z|$ as required. $\Box$

\medskip

\noindent For later reference, we recall the Moore bound for irregular graphs.

\begin{theorem}[Alon, Hoory, and Linial \cite{ahl}]\label{theoremmoore}
If $G$ is a graph of order $n$, girth at least $2r$ for some integer $r$, and average degree $d$ at least $2$, then
$n\geq 2\sum\limits_{i=0}^{r-1}(d-1)^i$.
\end{theorem}
We also need the following numerical fact.

\begin{lemma}\label{lemma2}
For positive integers $p$ and $f$ with $p\geq 5$ and $2p-1\leq f\leq {p\choose 2}$, 
$$\left(1+\frac{2(f-p)}{f+p}\right)^{\lceil\frac{p}{2}\rceil+1}>f-p+1.$$
\end{lemma}
{\it Proof:} For $p\geq 17$,
it follows from $f\geq 2p-1$
that $1+\frac{2(f-p)}{f+p}\geq 1.64$, 
and, since $1.64^{\lceil\frac{p}{2}\rceil+1}>{p\choose 2}-p+1$,
the desired inequality follows for these values of $p$.
For the finitely many pairs $(p,f)$ with $5\leq p \leq 16$ and $2p-1\leq f\leq {p\choose 2}$, we verified it using a computer. 
$\Box$

\medskip

\noindent We proceed to the proof of (\ref{eb3}).

\begin{theorem}
If $G$ is a graph of girth $g$ at least $3$ and minimum degree $\delta$ at least $2$,
then $Z(G)\geq (g-2)(\delta-2)+2$.
\end{theorem}
{\it Proof:} For $g=3$, the inequality simplifies to the known fact $Z(G)\geq \delta(G)$,
and, for $g=4$, it has been shown in \cite{gprs}. Now, let $g\geq 5$.
Let $X$ be a set of $g-2$ vertices of $G$ with $|N_G(X)|=\delta_{g-2}(G)$,
and, let $N=N_G(X)$.
By the girth condition, the components of $G[X]$ are trees,
and no vertex in $N$ has more than one neighbor in any component of $G[X]$.

Let $K_1,\ldots,K_p$ be the vertex sets of the components of $G[X]$.

If $p\geq 3$, and there are two vertices in $N$
that both have neighbors in two distinct components of $G[X]$,
then $G$ contains a cycle of order at most $2+|K_i|+|K_j|\leq 2+(g-2)-(m-2)<g$
which is a contradiction.
Similarly, if $p=2$, and there are three vertices $u$, $v$, and $w$ in $N$
that both have neighbors in $K_1$ and $K_2$,
then let $u_i$, $v_i$, and $w_i$ denote the corresponding neighbors in $K_i$ for $i\in [2]$, respectively.
If $G[K_1]$ contains a path between two of the vertices $u_1$, $v_1$, and $w_1$
avoiding the third, then $G$ contains a cycle of order at most $2+(|K_1|-1)+|K_2|=g-1$, which is a contradiction.
By symmetry, this implies $u_1=v_1=w_1$ and $u_2=v_2=w_2$,
and $G$ contains the cycle $u_1uu_2vu_1$ of order $4$, which is a contradiction.

Combining these observations, we obtain
\begin{eqnarray}\label{e1}
\sum\limits_{\{ i, j\}\in {[p]\choose 2}}|N_G(K_i)\cap N_G(K_j)|
\leq 
\begin{cases}
{p\choose 2} & \mbox{, for $p\geq 3$, and }\\
2p-2 & \mbox{, for $p\in \{ 1,2\}$}.
\end{cases}
\end{eqnarray}
Let the bipartite graph $H$ arise from $G[X\cup N]$ 
by contracting each component of $G[X]$ to a single vertex, 
and removing all edges of $G[N]$.
Lemma \ref{lemma1} and a simple counting implies
\begin{eqnarray*}
Z(G) & \geq & \delta_{g-2}(G)\\ 
& = & |N|\\
& = & \sum_{u\in V(H)\setminus N}d_H(u)-\sum_{v\in N}(d_H(v)-1)\\
& \geq & \sum_{i=1}^p\Big(\delta|K_i|-2(|K_i|-1)\Big)-\sum_{v\in N}(d_H(v)-1)\\
& = & (g-2)(\delta-2)+2+\left((2p-2)-\sum_{v\in N}(d_H(v)-1)\right).
\end{eqnarray*}
In view of (\ref{eb3}), we may assume $f\geq 2p-1$ for $f=\sum\limits_{v\in N}(d_H(v)-1)$.

Since 
$$2p-1\leq f=\sum\limits_{v\in N}(d_H(v)-1)
\leq \sum\limits_{v\in N}{d_H(v)\choose 2}
=\sum\limits_{\{ i, j\}\in {[p]\choose 2}}|N_G(K_i)\cap N_G(K_j)|,$$
(\ref{e1}) implies $p\geq 5$.

Let $H'$ arise by removing all vertices of degree $1$ from $H$.
Since every vertex $u$ in $V(H)\setminus N$ satisfies
$d_H(u)\geq \delta |K_i|-2(|K_i|-1)\geq 2$ for some $i\in [p]$,
the graph $H'$ contains all $p$ vertices of $V(H)\setminus N$.
Let $H'$ contain $q$ vertices of $N$.
Since $H'$ has order $p+q$ and size
$$\sum\limits_{v\in N\cap V(H')}d_H(v)
=q+\sum\limits_{v\in N}(d_H(v)-1)=q+f,$$
its average degree is at least $\frac{2(f+q)}{p+q}$,
which is at least $2$, because $f\geq 2p-1\geq p$.

If $H'$ contains a cycle of order $2\ell$,
then $G$ contains a cycle of order at most $(g-2)-(p-\ell)+\ell$.
By the girth condition, 
this implies that the bipartite graph $H'$ has girth at least 
$p+2$, if $p$ is even, and $p+3$, if $p$ is odd.

Using Theorem \ref{theoremmoore} and $f\geq q$, we obtain
\begin{eqnarray*}
p+q & \geq & 2\sum\limits_{i=0}^{\lceil\frac{p}{2}\rceil}\left(\frac{2(f+q)}{p+q}-1\right)^i\\
&=& 2\frac{p+q}{2f-2p}\left(\left(1+\frac{2(f-p)}{p+q}\right)^{\lceil\frac{p}{2}\rceil+1}-1\right)\\
&\geq &2\frac{p+q}{2f-2p}\left(\left(1+\frac{2(f-p)}{p+f}\right)^{\lceil\frac{p}{2}\rceil+1}-1\right),
\end{eqnarray*}
which implies
$\left(1+\frac{2(f-p)}{f+p}\right)^{\lceil\frac{p}{2}\rceil+1}\leq f-p+1$.
Since $f\geq 2p-1$, and, by (\ref{e1}), $f\leq {p\choose 2}$,
this contradicts Lemma \ref{lemma2},
which completes the proof.
$\Box$

\end{document}